%% file: 1kostant.tex
\documentclass[a4paper,12pt,reqno]{amsart}

\usepackage{graphicx}
\usepackage{amsmath}
\usepackage{amssymb}
\usepackage{amsfonts}
\usepackage{epsfig}

\makeatletter

 \@addtoreset{equation}{section}
 \@addtoreset{figure}{section}
 \@addtoreset{table}{section}
 \@addtoreset{theorem}{subsection}

\newtheorem{theorem}{Theorem}[section]

\newtheorem{corollary}[theorem]{Corollary}
\newtheorem{remark}[theorem]{Remark}
\newtheorem{proposition}[theorem]{Proposition}

\def\PerfProof{{\it Proof.\ }}

\makeatother

\makeindex

\begin{document}

\title[Kostant's generating functions, Ebeling's theorem
       and McKay's observation]{Kostant's generating functions, Ebeling's theorem
       and McKay's observation relating the Poincar\'{e} series}
       \author{Rafael Stekolshchik}
         \thanks{{\large email: rs2@biu.013.net.il}}

\date{}

\begin{abstract}
 We generalize B.~Kostant's construction of generating functions
to the case of multiply-laced diagrams and we prove for this case
W.~Ebeling's theorem which connects the Poincar\'{e} series
$[P_G(t)]_0$ and the Coxeter transformations. According to
W.~Ebeling's theorem
\begin{equation*}
        [P_G(t)]_0 = \frac{\mathcal{X}(t^2)}{\tilde{\mathcal{X}}(t^2)},
\end{equation*}
where $\mathcal{X}$ is the characteristic polynomial of the Coxeter
transformation and $\tilde{\mathcal{X}}$ is the characteristic
polynomial of the corresponding affine Coxeter transformation.

 We prove McKay's observation relating the Poincar\'{e} series
 $[P_G(t)]_i$:
 \begin{equation*}
    (t+t^{-1})[P_G(t)]_i = \sum\limits_{i \leftarrow j}[P_G(t)]_j,
 \end{equation*}
where $j$ runs over all vertices adjacent to $i$.
\end{abstract}

\vspace{7mm} 


\maketitle


\tableofcontents
\newpage


\input 2mckay.tex


\renewcommand{\appendixname}{}

  \printindex

\end{document}

%% file: 2mckay.tex
\section{Introduction}

B.~Kostant's construction of a vector-valued generating function
$P_G(t)$  appears in the context of the McKay correspondence
\cite{McK80} and gives a way to obtain multiplicities of
indecomposable representations $\rho_i$ of the binary polyhedral
group $G$ in the decomposition of $\pi_n|G$, \cite{Kos84}.

 Using B.~Kostant's construction we generalize
to the case of multiply-laced diagrams W.~Ebeling's theorem which
connects the Poincar\'{e} series $[P_G(t)]_0$ and the Coxeter
transformations. According to W.~Ebeling's theorem \cite{Ebl02}
\begin{equation*}
        [P_G(t)]_0 = \frac{\mathcal{X}(t^2)}{\tilde{\mathcal{X}}(t^2)},
\end{equation*}
where $\mathcal{X}$ is the characteristic polynomial of the Coxeter
transformation and $\tilde{\mathcal{X}}$ is the characteristic
polynomial of the corresponding affine Coxeter transformation.

 We prove McKay's observation \cite{McK99}
relating the Poincar\'{e} series $[P_G(t)]_i$:
 \begin{equation*}
    (t+t^{-1})[P_G(t)]_i = \sum\limits_{i \leftarrow j}[P_G(t)]_j,
 \end{equation*}
where $j$ runs over all vertices adjacent to $i$.

\subsection{The McKay correspondence}

Let $G$ be a finite subgroup of $SU(2)$, $\{\rho_0,\rho_1,\dots,
\rho_n\}$ be the set of all distinct irreducible finite dimensional
complex representations of $G$, of which $\rho_0$ is the trivial
one. Let $\rho: G \longrightarrow SU(2)$ be a faithful
 representation. Then, for each group $G$, we can define
 a matrix $A(G) = (a_{ij})$, by decomposing the tensor products:

\begin{equation}
 \label{main_McKay}
 \rho \otimes \rho_j = \bigoplus\limits_{k=0}^r a_{jk}\rho_k,
  \hspace{7mm} j = 0,1,...,r,
\end{equation}
where $a_{jk}$ is the multiplicity of $\rho_k$ in $\rho \otimes
\rho_j$. McKay \cite{McK80} observed that
\begin{center}
{\it The matrix $2I - A(G)$ is the Cartan matrix of the extended
  Dynkin diagram $\tilde\varDelta(G)$ associated to $G$.
  There is one-to-one correspondence between finite subgroups
  of $SU(2)$ and simply-laced extended Dynkin diagrams.
 }\end{center}

This remarkable observation, called the {\it McKay correspondence},
was based first on an explicit verification \cite{McK80}.

For the multiply-laced case, the McKay correspondence was extended
by D.~Happel, U.~Preiser, and C.~M.~Ringel in \cite{HPR80}, and by
P.~Slodowy in \cite[App.III]{Sl80}. We consider P.~Slodowy's
approach in \S\ref{slodowy}.

The systematic proof of the McKay correspondence based on the study
of affine Coxeter transformations was given by R.~Steinberg in
\cite{Stb85}.

 Other proofs of the McKay correspondence were given by
 G.~Gonzalez-Sprinberg and J.-L.~Verdier in \cite{GV83}, by
H.~Kn\"{o}rrer in \cite{Kn85}. A nice review is given by J.~van
Hoboken in \cite{Hob02}.

\subsection{The Slodowy generalization of the McKay correspondence}
 \label{slodowy}

\label{section_slodowy} Let us fix a pair $H \triangleleft G$ from
Table \ref{group_pairs}. We formulate now the essence of the Slodowy
correspondence \cite[App.III]{Sl80}.

\begin{table} [h]
  \centering
  \vspace{2mm}
  \caption{\hspace{3mm}The pairs $H \triangleleft G$ of binary polyhedral groups}
  \renewcommand{\arraystretch}{1.5} 
 \begin{tabular} {||c|c|c|c|c||}
 \hline \hline
   \quad Subgroup \quad &
   \quad Dynkin  \quad  &
   \quad Group \quad    &
   \quad Dynkin \quad   &
   \quad Index \quad    \cr
   \quad $H$  \quad     &
   \quad diagram $\varGamma(H)$ \quad &
   \quad $R$  \quad &
   \quad diagram $\varGamma(G)$ \quad &
   \quad $[G:H]$ \quad \\
 \hline \hline
      & & & & \cr
      $\mathcal{D}_2$ & ${D}_4$ & $\mathcal{T}$ & ${E}_6$ & 3 \cr
      & & & & \\
  \hline
      & & & & \cr
      $\mathcal{T}$  & ${E}_6$ & $\mathcal{O}$ & ${E}_7$ & 2 \cr
      & & & & \\
  \hline
      & & & & \cr
      $\mathcal{D}_{r-1}$ & ${D}_{n+1}$ & $ \mathcal{D}_{2(r-1)}$ &
      ${D}_{2n}$ & 2 \cr
      & & & & \\
 \hline
      & & & & \cr
      $\mathbb{Z}/2r\mathbb{Z}$ & ${A}_{n-1}$ & $\mathcal{D}_r$ &
      ${D}_{r+2}$ & 2 \cr
      & & & & \\
 \hline \hline
 \end{tabular}
  \label{group_pairs}
\end{table}

 1) Let $\rho_i$, where $i = 1,\dots,n$, be
irreducible representations of $G$; let $\rho^\downarrow_i$ be the
corresponding restricted representations of the subgroup $H$. Let
$\rho$ be a faithful representation of $H$, which may be considered
as the restriction of the fixed faithful representation $\rho_{f}$
of $G$. Then the following decomposition formula makes sense
\begin{equation}
 \label{lab_slodowy_1}
     \rho \otimes \rho^\downarrow_i =
          \bigoplus a_{ji} \rho^\downarrow_j
\end{equation}
and uniquely determines an $n\times{n}$ matrix $\tilde{A} =
(a_{ij})$ such that
\begin{equation}
   K = 2I - \tilde{A}
\end{equation}
(see \cite[p.163]{Sl80}), where $K$ is the Cartan matrix of the
corresponding folded extended Dynkin diagram given in Table
\ref{pairs_and_folded_diagr}

2) Let $\tau_i$, where $i = 1,\dots,n$, be irreducible
representations of the subgroup $H$, let $\tau^\uparrow_i$ be the
induced representations of the group $G$. Then the following
decomposition formula makes sense
\begin{equation}
 \label{lab_slodowy_2}
     \rho \otimes \tau^\uparrow_i =
          \bigoplus a_{ij} \tau^\uparrow_j,
\end{equation}
i.e., the decomposition of the induced representation is described
by the matrix $A^\vee = A^t$ which satisfies the relation
\begin{equation}
   K^\vee = 2I - \tilde{A}^\vee
\end{equation}
(see \cite[p.164]{Sl80}), where $K^\vee$ is the Cartan matrix of the
dual folded extended Dynkin diagram given in Table
\ref{pairs_and_folded_diagr}.

 We call matrices $\tilde{A}$
and $\tilde{A}^\vee$ the {\it Slodowy matrices}, they are analogs of
the McKay matrix. The {\it Slodowy correspondence} is an analogue to
the McKay correspondence for the multiply-laced case, so one can
speak about the {\it McKay-Slodowy correspondence}.

\begin{table} [h]
  \centering
  \renewcommand{\arraystretch}{1.5} 
  \vspace{2mm}
  \caption{\hspace{3mm}The pairs $H \triangleleft G$ and folded extended Dynkin
diagrams}
 \begin{tabular} {||c|c|c||}
 \hline \hline
   \quad Groups \quad  &
   \quad Dynkin diagram \quad  &
   \quad Folded extended  \quad \cr
   \quad $H \triangleleft G$ \quad &
   \quad $\tilde\varGamma(H)$ and $\tilde\varGamma(G)$ \quad &
   \quad Dynkin diagram \quad \\
 \hline \hline
      & & \cr
      $\mathcal{D}_2 \triangleleft \mathcal{T}$ &
      ${D}_4$ and ${E}_6$ &
      $\tilde{G}_{21}$ and $\tilde{G}_{22}$ \cr
      & & \\
 \hline
      & & \cr
      $\mathcal{T} \triangleleft \mathcal{O}$ &
      ${E}_6$ and ${E}_7$ &
      $\tilde{F}_{41}$ and $\tilde{G}_{42}$ \cr
      & & \\
 \hline
      & & \cr
      $\mathcal{D}_{r-1} \triangleleft \mathcal{D}_{2(r-1)}$ &
      ${D}_{n+1}$ and ${D}_{2n}$ &
      $\widetilde{DD}_{n}$ and $\widetilde{CD}_{n}$  \cr
      & & \\
 \hline
      & & \cr
      $\mathbb{Z}/2r\mathbb{Z} \triangleleft \mathcal{D}_r$ &
      ${A}_{n-1}$ and ${D}_{r+2}$ &
      $\tilde{B}_{n}$ and $\tilde{C}_{n}$  \cr
      & & \\
 \hline \hline
 \end{tabular}
  \label{pairs_and_folded_diagr}
\end{table}

\subsection{The Kostant generating function and Poincar\'{e} series}
\label{generating_fun}

Let ${\rm Sym}(\mathbb{C}^2)$ be the symmetric algebra over
$\mathbb{C}^2$, in other words, ${\rm Sym}(\mathbb{C}^2) =
\mathbb{C}[x_1, x_2]$. The symmetric algebra ${\rm
Sym}(\mathbb{C}^2)$ is a graded $\mathbb{C}$-algebra, see
\cite{Sp77}, \cite{Ben93}:
\begin{equation}
    {\rm Sym}(\mathbb{C}^2) =
\mathop{\oplus}\limits_{n=0}^{\infty}{\rm Sym}^n(\mathbb{C}^2).
\end{equation}

Let $\pi_n$ be the representation of $SU(2)$ in ${\rm
Sym}^n(\mathbb{C}^2)$ induced by its action on $\mathbb{C}^2$. The
set $\{\pi_n\}$, where $n = 0,1,...$ is the set of all irreducible
representations of $SU(2)$, see, e.g., \cite[\S37]{Zhe73}. Let $G$
be any finite subgroup of $SU(2)$. B.~Kostant in \cite{Kos84}
considered the question:

\medskip
{\sl how does $\pi_n | G$ decompose for any $n \in \mathbb{N}$?}
\medskip

The answer --- the decomposition $\pi_n | G$ --- is as follows:
\begin{equation}
  \label{decomp_pi_n}
   \pi_n | G = \sum\limits_{i=0}^r{m_i(n)\rho_i},
\end{equation}
where $\rho_i$ are irreducible representations of $G$, considered in
the context of {\it McKay correspondence}, see \cite{Kos84}). Thus,
the decomposition (\ref{decomp_pi_n}) reduces the question to the
following one:

\medskip
{\sl what are the multiplicities $m_i(n)$ equal to?}
\medskip

B.~Kostant in \cite{Kos84} obtained the multiplicities $m_i(n)$ by
means the orbit structure of the Coxeter transformation on the
highest root of the corresponding Lie algebra. For further details
concerning this orbit structure and the multiplicities $m_i(n)$, see
\S\ref{orbit_str}.

Note, that multiplicities $m_i(n)$ in (\ref{decomp_pi_n}) are
calculated as follows:
\begin{equation}
  \label{multipl_n}
     m_i(n) = <\pi_n | G , \rho_i>,
\end{equation}
(for the definition of the inner product $<\cdot , \cdot>$, see
(\ref{inner_prod}) ).
\begin{remark}
  \label{remark_decomp}
 {\rm
For further considerations, we extend the relation for
multiplicity (\ref{multipl_n}) to the cases of {\it restricted
representations}
  $\rho_i^\downarrow: = \rho_i\downarrow_H^G$
  and {\it induced representations}
   $\rho_i^\uparrow: = \rho_i\uparrow_H^G$, where $H$ is any subgroup of $G$
  (see \cite[Section A.5.1]{St05}):

\begin{equation}
  \label{decomp_pi_n_2}
   m_i^\downarrow(n) = <\pi_n | H ,\rho_i^\downarrow>,   \hspace{3mm}
   m_i^\uparrow(n) = <\pi_n | G , \rho_i^\uparrow>.
\end{equation}

We do not have any decomposition like (\ref{decomp_pi_n}) neither
for restricted representations $\rho_i^\downarrow$ nor for induced
representations $\rho_i^\uparrow$. Nevertheless, we will sometimes
denote both multiplicities $m_i^\downarrow(n)$ and
$m_i^\uparrow(n)$ in (\ref{decomp_pi_n_2}) by $m_i(n)$ as in
(\ref{decomp_pi_n}). }
\end{remark}

\begin{remark} {\rm
 \label{triv_repr}
1) A representation $\rho: G \longrightarrow GL_k(V)$ defines a
$k$-linear action $G$ on $V$ by
\begin{equation}
  gv = \rho(g)v.
\end{equation}
The pair $(V,\rho)$ is called a $G$-{\it module}. The case where
$\rho(g) = Id_V$ is called the {\it trivial representation} in
$V$. In this case
\begin{equation}
 \label{trivial}
  gv = v \text{ for all } g \in V.
\end{equation}
In (\ref{decomp_pi_n}), the trivial representation $\rho_0$
corresponds to a particular vertex (see \cite{McK80}),
which extends the Dynkin diagram to the extended Dynkin diagram.

2) Let $\rho_0(H)$ (resp. $\rho_0(G)$) be the trivial
representation of any subgroup $H \subset G$ (resp. of group $G$).
The trivial representation $\rho_0(H): H \longrightarrow GL_k(V)$
coincides with the {\it restricted representation}
$\rho_0\downarrow_H^G: G \longrightarrow GL_k(V)$, and the trivial
representation $\rho_0(G): G \longrightarrow GL_k(V)$ coincides
with the {\it induced representation} $\rho_0\uparrow_H^G: H
\longrightarrow GL_k(V)$.}
\end{remark}

Since there is one-to-one correspondence between the $\rho_i$ and
the vertices of the Dynkin diagram, we can define (see
\cite[p.211]{Kos84}) the vectors $v_n$, where $n\in\mathbb{Z}_+$,
as follows:
\begin{equation}
 \label{def_vn}
 v_n  = \sum\limits_{i=0}^r{m_i(n)}\alpha_i,
     \text{ where } \pi_n | G = \sum\limits_{i=0}^r{m_i(n)\rho_i}, 
\end{equation}
where $\alpha_i$ are simple roots of the corresponding extended
Dynkin diagram. Similarly, for  the multiply-laced case, we define
vectors $v_n$ to be:
\begin{equation}
 \label{def_vn_1}
   v_n  = \sum\limits_{i=0}^r{m_i^\uparrow(n)}\alpha_i \hspace{5mm} \text{
   or } \hspace{5mm}
  v_n  = \sum\limits_{i=0}^r{m_i^\downarrow(n)}\alpha_i,
\end{equation}
where the multiplicities $m_i^\uparrow(n)$ and $m_i^\downarrow(n)$
are defined by (\ref{decomp_pi_n_2}). The vector $v_n$ belongs to
the root lattice generated by simple roots. Following B.~Kostant,
we define the generating function $P_G(t)$ for cases
(\ref{def_vn}) and (\ref{def_vn_1}) as follows:
\begin{equation}
 \label{Kostant_gen_func}
    P_G(t) = ([P_G(t)]_0, [P_G(t)]_1, \ldots , [P_G(t)]_r)^t
           := \sum\limits_{n=0}^{\infty}v_n{t^n},
\end{equation}
the components of the vector $P_G(t)$ being the following series
\begin{equation}
 \label{gen_func_i}
    [P_G(t)]_i = \sum\limits_{n=0}^{\infty}\tilde m_i(n){t^n},
\end{equation}
where $i = 0,1,\dots,r$ and $\tilde m_i(n)$ designates $m_i(n),
m_i^\uparrow(n)$ or $m_i^\downarrow(n)$.  In particular, for $i=0$,
we have
\begin{equation}
    [P_G(t)]_0 = \sum\limits_{n=0}^{\infty}m_0(n){t^n},
\end{equation}
where $m_0(n)$ is the multiplicity of the trivial representation
$\rho_0$ (Remark \ref{triv_repr}) in ${\rm Sym}^n(\mathbb{C}^2)$.
The {\it algebra of invariants} $R^G$ is a subalgebra of the {\it
symmetric algebra} ${\rm Sym}(\mathbb{C}^2)$.
 Thanks to (\ref{trivial}), we see that
$R^G$ coincides with ${\rm Sym}(\mathbb{C}^2)$, and $[P_G(t)]_0$ is
the Poincar\'{e} series of the algebra of invariants ${\rm
Sym}(\mathbb{C}^2)^G$, i.e.,
\begin{equation}
  \label{poincare_alg_inv}
    [P_G(t)]_0 = P({\rm Sym}(\mathbb{C}^2)^G,t).
\end{equation}
(see \cite[p.221, Rem.3.2]{Kos84}).
\begin{remark}
\label{poincare_m_case} {\rm According to Remark \ref{triv_repr},
heading 2), we have
\begin{equation}
  \label{poincare_alg_inv_2}
    [P_H(t)]_0 = P({\rm Sym}(\mathbb{C}^2)^{\rho_0\downarrow_H^G},t), \hspace{7mm}
    [P_G(t)]_0 = P({\rm Sym}(\mathbb{C}^2)^{\rho_0\uparrow_H^G},t).
\end{equation}
We will need this fact in the proof of W.~Ebeling's theorem for the
multiply-laced case in \S\ref{ebeling}, see Remark
 \ref{remark_decomp}.
 }
\end{remark}
The following theorem gives a remarkable formula for calculating the
Poincar\'{e} series for the binary polyhedral groups. The theorem is
known in different forms, see Kostant \cite{Kos84}, Kn\"{o}rrer
\cite{Kn85}, Gonsales-Sprinberg, Verdier \cite{GV83}. B.~Kostant in
\cite{Kos84} shows it in the context of the Coxeter number $h$.

\begin{theorem} 
  The Poincar\'{e} series $[P_G(t)]_0$ can be calculated as the
following rational function:
\begin{equation}
  \label{K_K_GV}
      [P_G(t)]_0 = \frac{1 + t^h}{(1 - t^a)(1 - t^b)},
\end{equation}
where
\begin{equation}
  \label{Kostant_numbers_a_b}
      b = h + 2 - a, \text{ and } ab = 2|G|.
\end{equation}
\end{theorem}

For a proof, see Theorem 1.4 and Theorem 1.8 from \cite{Kos84},
\cite[p.185]{Kn85}, \cite[p.428]{GV83}. We call the numbers $a$ and
$b$ the {\it Kostant numbers}. They can be easily calculated, see
Table \ref{Kostant_numbers}, compare also with \cite[Table
A.2]{St05}. Note, that $a = 2d$, where $d$ is the maximal coordinate
of the {\it nil-root} vector from the kernel of the Tits form,
\cite{Kac80}.

\begin{table} [h]
  \centering
  \vspace{2mm}
  \caption{\hspace{3mm}The binary polyhedral groups (BPG) and the Kostant numbers
$a$, $b$}
   \renewcommand{\arraystretch}{1.5} 
  \begin{tabular} {||c|c|c|c|c|c||}
  \hline \hline
     \quad  Dynkin \quad    &
     \quad Order of \quad   &
     \quad BPG \quad        &
     \quad Coxeter \quad    &
     \quad $a$ \quad        &
     \quad $b$ \quad        \cr
      diagram      &  group   &            &  number        &    &   \\
  \hline \hline
      & & & & & \cr
      ${A}_{n-1}$  & $n$     & $\mathbb{Z}/n\mathbb{Z}$ & $n$ &  2 & $n$ \cr
      & & & & & \\
  \hline
      & & & & & \cr
      ${D}_{n+2}$  & $4n$    & $\mathcal{D}_n$ &   $2n+2$   & 4  & $2n$ \cr
      & & & & & \\
  \hline
      & & & & & \cr
      ${E}_6$  & 24  & $\mathcal{T}$   &   12           & 6  & 8 \cr
      & & & & & \\
  \hline
      & & & & & \cr
      ${E}_7$  & 48  & $\mathcal{O}$   &   18           & 8  & 12 \cr
      & & & & & \\
  \hline
      & & & & & \cr
      ${E}_8$  & 120 & $\mathcal{J}$   &   30           & 12 & 20 \cr
      & & & & & \\
  \hline  \hline
  \end{tabular}
    \label{Kostant_numbers}
\end{table}

\subsection{The characters and the McKay operator}
 \label{sect_char_McKay}

  Let $\chi_1, \chi_2,\dots \chi_r$ be all irreducible $\mathbb{C}$-characters
of a finite group $G$ corresponding to irreducible representations
$\rho_1, \rho_2, \dots, \rho_r$, and let $\chi_1$ correspond to
the trivial representation, i.e., $\chi_1(g) = 1$ for all $g \in
G$.

 All characters constitute the {\it
character algebra} $C(G)$ of $G$ since $C(G)$ is also a vector
space  over $\mathbb{C}$. An {\it hermitian inner product} $<\cdot
, \cdot>$ on $C(G)$ is defined as follows. For characters $\alpha,
\beta \in C(G)$, let
\begin{equation}
  \label{inner_prod}
   <\alpha, \beta> ~=~ \frac{1}{|G|}
     \sum\limits_{g \in G}\alpha(g)\overline{\beta(g)}
\end{equation}
Sometimes, we will write {\it inner product} $<\rho_i, \rho_j>$ of
the representations meaning actually the inner product of the
corresponding characters $<\chi_{\rho_i}, \chi_{\rho_j}>$.
 Let $z_{ijk} = <\chi_i\chi_j, \chi_k>$, where
 $\chi_i\chi_j$ corresponds to the representation $\rho_i\otimes\rho_j$.
It is known that
$z_{ijk}$ is the multiplicity of the representation $\rho_k$  in
$\rho_i\otimes\rho_j$ and $z_{ijk} = z_{jik}$.
The numbers $z_{ijk}$ are integer and are called the {\it structure constants},
see, e.g., \cite[p.765]{Kar92}.

For every $i \in \{1,...,r\}$, there exists some $\stackrel{\wedge}i \in
\{1,...,r\}$
such that
\begin{equation}
  \chi_{\stackrel{\wedge}{i}}(g) = \overline{\chi_i(g)}
  \text{ for all } g \in G.
\end{equation}

 The character
$\chi_{\stackrel{\wedge}{i}}$ corresponds to the {\it
contragredient representation} $\rho_{\stackrel{\wedge}{i}}$
determined from the relation
\begin{equation}
  \label{def_contragr}
    \rho_{\stackrel{\wedge}{i}}(g) = \rho_i(g)^{'{-1}}.
\end{equation}
We have
\begin{equation}
  \label{char_contragr}
   <\chi_i\chi_j, \chi_k> =
   <\chi_i, \chi_{\stackrel{\wedge}{j}}\chi_k>
\end{equation}
since
\begin{equation*}
 \begin{split}
  <\chi_i\chi_j, \chi_k> = &
   \frac{1}{|G|}\sum\limits_{g \in G}
  \chi_i(g)\chi_j(g)\overline{\chi_k(g)} = \\
  & \frac{1}{|G|}\sum\limits_{g \in G}
  \chi_i(g)(\overline{\overline{\chi_j(g)}\chi_k(g)}) =
   <\chi_i, \chi_{\stackrel{\wedge}{j}}\chi_k>.
 \end{split}
\end{equation*}

\begin{remark}
\label{group_su2}
{\rm
The group $SU(2)$ is the set of all unitary unimodular $2\times2$ matrices $u$,
  i.e.,
\begin{equation*}
  \label{unit_unim}
  \begin{split}
    &  u^* = u^{-1}  \text{ ({\it unitarity}) }, \\
    &  \det(u) = 1   \text{ ({\it unimodularity}) }.
  \end{split}
\end{equation*}
The matrices $u \in SU(2)$ have the following form:
\begin{equation}
  \label{su2_matr}
  u = \left (
       \begin{array}{cc}
         a     &  b \\
         -b^*  &  a^* \\
       \end{array}
      \right ), \text{ and }
   u^* =      \left (
       \begin{array}{cc}
         a^*  &  -b \\
         b^*  &  a \\
       \end{array}
      \right ),
      \text{ where } aa^* + bb^* = 1,
\end{equation}
see, e.g., \cite[Ch.9, \S6]{Ha89}. The mutually inverse matrices
  $u$ and $u^{-1}$ are
\begin{equation}
  \label{sl2_matr}
  u = \left (
       \begin{array}{cc}
         a  &  b \\
         c  &  d \\
       \end{array}
      \right ), \text{ and }
  u^{-1}=  \left (
       \begin{array}{cc}
         d  &  -b \\
         -c  &  a \\
       \end{array}
      \right ), \text{ where } ad - bc = 1.
\end{equation}

Set
\begin{equation}
  s = \left (
       \begin{array}{cc}
         0  &  -1 \\
         1  &  0 \\
       \end{array}
      \right ), \text{ then }
  s^{-1} = s^3 = \left (
       \begin{array}{cc}
         0  &  1 \\
        -1  &  0 \\
       \end{array}
      \right ).
\end{equation}
For  any $u \in SU(2)$, we have
\begin{equation}
 \label{elem_weyl}
  sus^{-1} = u^{~'{-1}}.
\end{equation}
The element $s$ is called the {\it Weyl element}. }
\end{remark}

According to (\ref{def_contragr}) and (\ref{elem_weyl}) we see
that every finite dimensional representation of the group $SL(2,
\mathbb{C})$ (and hence, of $SU(2)$) is equivalent to its
contragredient representation, see \cite[\S37, Rem.3]{Zhe73}. Thus
by (\ref{char_contragr}), for representations $\rho_i$ of any
finite subgroup $G \subset SU(2)$,  we have
\begin{equation}
  \label{char_contragr_1}
   <\chi_i\chi_j, \chi_k> ~=~ <\chi_i, \chi_j\chi_k>.
\end{equation}

 Relation (\ref{char_contragr_1}) holds
also for characters of restricted representations $\chi^\downarrow
:= \chi\downarrow^G_H$ and induced representations
 $\chi^\uparrow := \chi\uparrow^G_H$, see Remark \ref{remark_decomp}:

\begin{equation}
  \label{char_contragr_2}
\renewcommand{\arraystretch}{1.2}
\begin{array}{l}
<\chi_i\chi^\downarrow_j, \chi_k>_H ~=~ <\chi_i, \chi^\downarrow_j\chi_k>_H,\\
   <\chi_i\chi^\uparrow_j, \chi_k>_G ~=~ <\chi_i, \chi^\uparrow_j\chi_k>_G.
\end{array}
\end{equation}
Indeed, every restricted representation $\chi^\downarrow_j$ (resp.
induced representation $\chi^\uparrow_j$) is decomposed into the
direct sum of irreducible characters $\chi_s \in Irr(H)$ (here,
$Irr(H)$ is the set of irreducible characters of $H$) with some
integer coefficients $a_s$, for example:
\begin{equation}
  \label{char_contragr_3}
   \chi^\downarrow_j = \sum\limits_{\chi_s \in Irr(H)}{a_s}\chi_s,
\end{equation}
and since $\overline{a}_s = a_s$ for all $\chi_s \in Irr(H)$, we
have
\begin{equation}
\begin{split}
  \label{char_contragr_4}
   <\chi_i\chi^\downarrow_j, \chi_k>_H = &
    \sum\limits_{\chi_s \in Irr(H)} <a_s\chi_i\chi_s, \chi_k>_H = \\
   & \sum\limits_{\chi_s \in Irr(H)} a_s<\chi_i, \chi_s\chi_k>_H =  \\
   &  <\chi_i, \sum\limits_{\chi_s \in Irr(H)}a_s\chi_s\chi_k>_H =
    <\chi_i, \chi^\downarrow_j\chi_k>_H.
\end{split}
\end{equation}

 The matrix of multiplicities
$A := A(G)$ from (\ref{main_McKay}) was introduced by J.~McKay in
\cite{McK80}; it plays the central role in the {\it McKay
correspondence}, see \cite{McK80}. We call this matrix --- or the
corresponding operator
--- the {\it McKay matrix} or the {\it McKay operator}.

Similarly, let $\tilde{A}$ and $\tilde{A}^\vee$ be matrices
 of multiplicities (\ref{lab_slodowy_1}),  (\ref{lab_slodowy_2}).
These matrices were introduced by P.~Slodowy \cite{Sl80} by analogy
with the McKay matrix for the multiply-laced case. We call these
matrices the {\it Slodowy operators}.

The following result of B.~Kostant \cite{Kos84}, which holds for
the McKay operator holds also for the Slodowy operators.

\begin{proposition}
  \label{kostant_prop}
  If $B$ is either the McKay operator $A$ or the Slodowy operator
  $\tilde{A}$ or $\tilde{A}^\vee$, then
 \begin{equation}
   \label{Kostant_relation}
      Bv_n = v_{n-1} + v_{n+1}.
 \end{equation}

\end{proposition}
\PerfProof
From now on
\begin{equation}
    \rho_i =
     \begin{cases}
        \rho_i &\text{ for } B = A, \vspace{2mm} \\
        \rho_i^\downarrow &\text{ for } B = \tilde{A}, \vspace{2mm} \\
        \rho_i^\uparrow &\text{ for } B = \tilde{A}^\vee, \\
     \end{cases}
 \hspace{7mm}
    m_i(n) =
     \begin{cases}
        m_i(n) &\text{ for } B = A, \vspace{2mm} \\
        m_i^\downarrow(n) &\text{ for } B = \tilde{A}, \vspace{2mm} \\
        m_i^\uparrow(n) &\text{ for } B = \tilde{A}^\vee. \\
     \end{cases}
\end{equation}

By (\ref{def_vn}), (\ref{def_vn_1}), and by definition of the McKay
operator (\ref{main_McKay}) and by definition of the Slodowy
operator (\ref{lab_slodowy_1}), (\ref{lab_slodowy_2}), we have
\begin{equation}
  \label{McKay_oper}
      Bv_n = B
    \left (
    \begin{array}{c}
       m_0(n) \\
       \dots  \\
       m_r(n)
    \end{array}
    \right ) =
    \left (
    \begin{array}{c}
       \sum a_{0i}m_i(n) \\
       \dots  \\
       \sum a_{ri}m_i(n) \\
    \end{array}
    \right ) =
    \left (
    \begin{array}{c}
       \sum a_{0i}<\rho_i, \pi_n> \\
       \dots  \\
       \sum a_{ri}<\rho_i, \pi_n> \\
    \end{array}
    \right ).
\end{equation}
By (\ref{main_McKay}), (\ref{lab_slodowy_1}) and
(\ref{lab_slodowy_2}) we have
\begin{equation*}
  \sum\limits_{i=1}^r a_{0i}<\rho_i, \pi_n> ~=~
  <\sum\limits_{i=1}^r a_{0i}\rho_i, \pi_n> ~=~
  <\rho_{reg}\otimes\rho_i, \pi_n>, \\
\end{equation*}
and from (\ref{McKay_oper}) we obtain

\begin{equation}
  \label{McKay_oper_1}
      Bv_n =
    \left (
    \begin{array}{c}
       <\rho_{reg}\otimes\rho_0, \pi_n> \\
       \dots  \\
       <\rho_{reg}\otimes\rho_r, \pi_n> \\
    \end{array}
    \right ).
\end{equation}
  Here $\rho_{reg}$ is the regular two-dimensional representation which
coincides with the representation $\pi_1$ in ${\rm
Sym}^2(\mathbb{C}^2)$ from \S\ref{generating_fun}. Thus,
\begin{equation}
  \label{McKay_oper_2}
   Bv_n =
    \left (
    \begin{array}{c}
       <\pi_1\otimes\rho_0, \pi_n> \\
       \dots  \\
       <\pi_1\otimes\rho_r, \pi_n> \\
    \end{array}
    \right ).
\end{equation}
From (\ref{char_contragr_1}), (\ref{char_contragr_4}) we obtain
\begin{equation}
  \label{McKay_oper_3}
   Bv_n =
    \left (
    \begin{array}{c}
       <\rho_0, \pi_1\otimes\pi_n> \\
       \dots  \\
       <\rho_r, \pi_1\otimes\pi_n> \\
    \end{array}
    \right ).
\end{equation}
By Clebsch-Gordan formula we have
\begin{equation}
  \label{McKay_oper_4}
   \pi_1\otimes\pi_n = \pi_{n-1} \oplus \pi_{n+1},
\end{equation}
where $\pi_{-1}$ is the zero representation, see
\cite[exs.3.2.4]{Sp77} or \cite[Ch.5, \S6,\S7]{Ha89}. From
(\ref{McKay_oper_3}) and (\ref{McKay_oper_4}) we have
(\ref{Kostant_relation}). \qedsymbol

For the following corollary, see \cite[p.222]{Kos84} and also
\cite[\S4.1]{Sp87}.
\begin{corollary}
Let $x = P_G(t)$ be given by (\ref{Kostant_gen_func}). Then
\begin{equation}
  \label{McKay_oper_5}
  tBx = (1 + t^2)x - v_0, \\
\end{equation}
  where $B$ is either the McKay operator $A$ or the Slodowy operators
  $\tilde{A}$, $\tilde{A}^\vee$.
\end{corollary}
\PerfProof
  From (\ref{Kostant_relation}) we obtain

\begin{equation}
\begin{split}
 & Bx = \sum\limits_{n=0}^{\infty}Bv_n{t^n} =
        \sum\limits_{n=0}^{\infty}(v_{n-1} + v_{n+1}){t^n} =
        \sum\limits_{n=0}^{\infty}v_{n-1}t^n +
        \sum\limits_{n=0}^{\infty}v_{n+1}t^n =  \\
 &  t\sum\limits_{n=1}^{\infty}v_{n-1}t^{n-1} +
    t^{-1}\sum\limits_{n=0}^{\infty}v_{n+1}t^{n+1} =
    tx + t^{-1}(\sum\limits_{n=0}^{\infty}v_n{t^n} - v_0) = \\
 &  tx + t^{-1}x - t^{-1}v_0.  \qed
\end{split}
\end{equation}

\section{The Poincar\'{e} series and W.~Ebeling's theorem}
\label{ebeling}

W.~Ebeling in \cite{Ebl02} makes use of the Kostant relation
(\ref{Kostant_relation}) and deduces a new remarkable fact about
the Poincar\'{e} series, a fact that shows that the Poincar\'{e}
series of a binary polyhedral group (see (\ref{K_K_GV})) is the
quotient of two polynomials: the characteristic polynomial of the
Coxeter transformation and the characteristic polynomial of the
corresponding affine Coxeter transformation, see
\cite[Th.2]{Ebl02}.

We show W.~Ebeling's theorem also for the multiply-laced case, see
Theorem \ref{theorem_ebeling}. The Poincar\'{e} series for the
multiply-laced case is defined by (\ref{poincare_alg_inv_2}).

\begin{table} 
  \centering
  \vspace{2mm}
  \caption{\hspace{3mm}The characteristic polynomials $\mathcal{X}$, $\tilde{\mathcal{X}}$ and
the Poincar\'{e} series}
  \renewcommand{\arraystretch}{1.5}  
  \begin{tabular} {|c|c|c|c|}
  \hline \hline
    Dynkin   & Coxeter
             & Affine Coxeter
             & Quotient \cr
    diagram  & transformation $\mathcal{X}$
             & transformation $\tilde{\mathcal{X}}$
             & $p(\lambda) = \displaystyle\frac{\mathcal{X}}{\tilde{\mathcal{X}}}$ \\
  \hline \hline
     & & & \cr
       ${D}_{4}$
     & $(\lambda + 1)(\lambda^{3} + 1)$
     & $(\lambda - 1)^2(\lambda + 1)^3$
     & $\displaystyle\frac{\lambda^3 + 1}{(\lambda^2 - 1)^2}$ \cr
     & & & \\
  \hline
     & & & \cr
       ${D}_{n+1}$
     & $(\lambda + 1)(\lambda^{n} + 1)$
     & $(\lambda^{n-1} -  1)(\lambda - 1)(\lambda + 1)^2$
     & $\displaystyle\frac{\lambda^{n} + 1}{(\lambda^{n-1} -  1)(\lambda^2 -
1)}$ \cr
     & & & \cr
  \hline
     & & & \cr
       ${E}_6$
     & $\displaystyle\frac{(\lambda^{6} + 1)}{(\lambda^2 + 1)}
       \frac{(\lambda^{3} - 1)}{(\lambda - 1)}$
     & $(\lambda^3 -  1)^2(\lambda + 1)$
     & $\displaystyle\frac{\lambda^6 + 1}{(\lambda^4 -  1)(\lambda^3 - 1)}$ \cr
     & & & \cr
  \hline
     & & & \cr
       ${E}_7$
     & $\displaystyle\frac{(\lambda + 1)(\lambda^{9} + 1)}{(\lambda^3 + 1)}$
     & $(\lambda^4 -  1)(\lambda^3 -  1)(\lambda + 1)$
     & $\displaystyle\frac{\lambda^9 + 1}{(\lambda^4 -  1)(\lambda^6 - 1)}$ \cr
     & & & \cr
  \hline
     & & & \cr
       ${E}_8$
     & $\displaystyle\frac{(\lambda^{15} + 1)(\lambda + 1)}{(\lambda^5 +
1)(\lambda^3  + 1)}$
     & $(\lambda^5 -  1)(\lambda^3 -  1)(\lambda + 1)$
     & $\displaystyle\frac{\lambda^{15} + 1}{(\lambda^{10} -  1)(\lambda^6 -
1)}$ \cr
     & & & \cr
  \hline
     & & & \cr
       ${B}_n$
     & $\lambda^{n} +  1$
     & $(\lambda^{n-1} -  1)(\lambda^2 -1)$
     & $\displaystyle\frac{\lambda^{n} + 1}{(\lambda^{n-1} -  1)(\lambda^2 -
1)}$ \cr
  \hline
     & & & \cr
       ${C}_n$
     & $\lambda^{n} +  1$
     & $(\lambda^n -  1)(\lambda - 1)$
     & $\displaystyle\frac{\lambda^{n} + 1}{(\lambda^n -  1)(\lambda - 1)}$ \cr
  \hline
     & & & \cr
       ${F}_4$
     & $\displaystyle\frac{\lambda^{6} + 1}{\lambda^2 + 1}$
     & $(\lambda^2 -  1)(\lambda^3 - 1)$
     & $\displaystyle\frac{\lambda^6 + 1}{(\lambda^4 -  1)(\lambda^3 - 1)}$ \\
  \hline
     & & & \cr
       ${G}_2$
     & $\displaystyle\frac{\lambda^{3} + 1}{\lambda + 1}$
     & $(\lambda-1)^2(\lambda  + 1)$
     & $\displaystyle\frac{\lambda^3 + 1}{(\lambda^2 - 1)^2}$ \\
  \hline\hline
     & & & \cr
       ${A}_n$
     & $\displaystyle\frac{\lambda^{n+1} - 1}{\lambda - 1}$
     & $(\lambda^{n - k + 1} - 1)(\lambda^k - 1)$
     & $\displaystyle\frac{\lambda^{n+1} - 1}
       {(\lambda - 1)(\lambda^{n - k + 1} - 1)(\lambda^k - 1)}$ \\
  \hline
     & & & \cr
       ${A}_{2n-1}$
     & $\displaystyle\frac{\lambda^{2n} - 1}{\lambda - 1}$
     & $(\lambda^n - 1)^2 \text{ for } k = n$
     & $\displaystyle\frac{\lambda^n + 1}
       {(\lambda^n - 1)(\lambda - 1)}$ \\
  \hline \hline
\end{tabular}
  \label{table_char_polyn_and_Poincare}
\end{table}

\begin{theorem} [generalized W.Ebeling's theorem \cite{Ebl02}]
 \label{theorem_ebeling}
  Let $G$ be a binary polyhedral group and $[P_G(t)]_0$ the Poincar\'{e}
  series (\ref{poincare_alg_inv}) of the algebra of invariants
${\rm Sym}(\mathbb{C}^2)^G$. Then
\begin{equation}
       [P_G(t)]_0 = \frac{\det{M_0}(t)}{\det{M}(t)},
\end{equation}
where
\begin{equation}
   \det{M}(t) = \det|t^{2}I - {\bf C}_a|, \hspace{5mm}
   \det{M_0}(t) = \det|t^{2}I - {\bf C}|,
\end{equation}
{\bf C} is the Coxeter transformation and ${\bf C}_a$ is the
corresponding affine Coxeter transformation.
\end{theorem}

\PerfProof
By (\ref{McKay_oper_5}) we have
\begin{equation}
   [(1 + t^2)I - tB]x = v_0,
\end{equation}
where $x$ is the vector $P_G(t)$ and by Cramer's rule the first coordinate
$P_G(t)$ is
\begin{equation}
       [P_G(t)]_0 = \frac{\det{M_0}(t)}{\det{M}(t)},
\end{equation}
where
\begin{equation}
 \label{M_def_affine}
       \det{M}(t) = \det\left((1 + t^2)I - tB\right),
\end{equation}
and $M_0(t)$ is the matrix obtained by replacing the first column
of $M(t)$ by $v_0 = (1,0,...,0)^t$. The vector $v_0$ corresponds
to the trivial representation $\pi_0$, and by the McKay
correspondence, $v_0$ corresponds to the particular vertex which
extends the Dynkin diagram to the extended Dynkin diagram, see
Remark \ref{triv_repr} and (\ref{def_vn}). Therefore, if
$\det{M}(t)$ corresponds to the affine Coxeter transformation, and
\begin{equation}
 \label{M_affine_cox}
       \det{M}(t) = \det|t^{2}I - {\bf C}_a|,
\end{equation}
then $\det{M_0}(t)$ corresponds to the Coxeter transformation, and
\begin{equation}
 \label{M_cox}
       \det{M_0}(t) = \det|t^{2}I - {\bf C}|.
\end{equation}
So, it suffices to prove (\ref{M_affine_cox}), i.e.,
\begin{equation}
 \label{M_affine_cox_2}
        \det[(1 + t^2)I - tB] = \det|t^{2}I - {\bf C}_a|.
\end{equation}
If $B$ is the McKay operator $A$ given by (\ref{main_McKay}), then
\begin{equation}
 \label{McKay_operator}
       B = 2I - K =
    \left ( \begin{array}{cc}
            2I &  0  \\
            0  &  2I
            \end{array}
    \right ) -
    \left ( \begin{array}{cc}
            2I   & 2D \\
            2D^t &  2I
            \end{array}
    \right ) =
    \left ( \begin{array}{cc}
            0     & -2D \\
            -2D^t &  0
            \end{array}
    \right ),
\end{equation}
where $K$ is a symmetric Cartan matrix \cite[(3.2)]{St05}. If $B$ is
the Slodowy operator $\tilde{A}$ or $\tilde{A}^\vee$ given by
(\ref{lab_slodowy_1}), (\ref{lab_slodowy_2}), then
\begin{equation}
 \label{Slodowy_operator}
       B = 2I - K =
    \left ( \begin{array}{cc}
            2I &  0  \\
            0  &  2I
            \end{array}
    \right ) -
    \left ( \begin{array}{cc}
            2I   & 2D \\
            2F &  2I
            \end{array}
    \right ) =
    \left ( \begin{array}{cc}
            0     & -2D \\
            -2F &  0
            \end{array}
    \right ),
\end{equation}
where $K$ is the symmetrizable but not symmetric Cartan matrix
\cite[(3.4)]{St05}. Thus, in the generic case
\begin{equation}
 \label{M_t}
       M(t) = (1+t^2)I - tB =
    \left ( \begin{array}{cc}
            1+t^2 &  2tD     \\
            2tF  &  1 + t^2
            \end{array}
     \right ).
\end{equation}
Assuming $t \neq 0$ we deduce from (\ref{M_t}) that
\begin{equation}
 \label{M_t_2}
 \begin{array}{cc}
       M(t)
    \left ( \begin{array}{c}
            x    \\
            y
            \end{array}
    \right ) = 0 & \Longleftrightarrow
    \left \{
     \begin{array}{c}
            (1 + t^2)x = -2tDy,  \\
            2tFx = -(1 + t^2)y.
            \end{array}
     \right .  \vspace{5mm} \\
    & \Longleftrightarrow
    \left \{
     \begin{array}{c}
            \displaystyle\frac{(1 + t^2)^2}{4t^2}x = FDy,   \vspace{3mm} \\
            \displaystyle\frac{(1 + t^2)^2}{4t^2}y = DFy.
            \end{array}
     \right .
  \end{array}
\end{equation}
Here we use the Jordan form theory of the Coxeter transformations
constructed in \cite[Section 3]{St05}. According to
\cite[3.14]{St05},  \cite[Propositions 3.3 and 3.9]{St05} and we see
that $t^2$ is an eigenvalue of the affine Coxeter transformation
${\bf C}_a$, i.e., (\ref{M_affine_cox_2}) together with
(\ref{M_affine_cox}) are proved. \qedsymbol

For the results of calculations using W.~Ebeling's theorem, see
Table \ref{table_char_polyn_and_Poincare}.

\begin{remark} {\rm
1) The characteristic polynomials $\mathcal{X}$ for the Coxeter
transformation and $\tilde{\mathcal{X}}$ for the affine Coxeter
transformation in Table \ref{table_char_polyn_and_Poincare} are
taken from \cite[Tables 1.1 and 1.2]{St05}. Pay attention to the
fact that the affine Dynkin diagram for ${B}_n$ is
$\widetilde{CD}_n$, (\cite[Tab.2]{Bo}),  and the affine Dynkin
diagram for ${C}_n$ is $\tilde{C}_n$, (\cite[Tab.3]{Bo}), see
\cite[Fig.~2.2]{St05}.

 2) The characteristic
polynomial $\mathcal{X}$ for the affine Coxeter transformation of
${A}_n$ depends on the {\it index of the conjugacy class} $k$ of the
Coxeter transformation, see \cite[4.5]{St05}. In the case of ${A}_n$
(for every $k = 1, 2,..., n$) the quotient $p(\lambda) =
\displaystyle\frac{\mathcal{X}}{\tilde{\mathcal{X}}}$ contains three
factors in the denominator, and its form is different from
(\ref{K_K_GV}), see Table \ref{table_char_polyn_and_Poincare}.

For the case ${A}_{2n-1}$ and $k = n$, we have
\begin{equation}
 \begin{split}
     p(\lambda) = & \frac{\lambda^{2n} - 1}
                 {(\lambda - 1)(\lambda^{2n-k} - 1)(\lambda^k - 1)} = \vspace{2mm} \\
              & \frac{\lambda^{2n} - 1}
                 {(\lambda - 1)(\lambda^n - 1)(\lambda^n - 1)} =
               \frac{\lambda^n + 1}
                 {(\lambda^n - 1)(\lambda - 1)}
 \end{split}
\end{equation}
and $p(\lambda)$ again is of the form (\ref{K_K_GV}),
see Table \ref{table_char_polyn_and_Poincare}.

3) The quotients $p(\lambda)$ coincide for the following pairs:
\begin{equation}
   \begin{array}{cc}
       {D}_4  \text{ and } {G}_2, \hspace{7mm}
      & {E}_6  \text{ and } {F}_4, \\
       {D}_{n+1} \text{ and } {B}_n  (n \geq 4), \hspace{7mm}
      &  {A}_{2n-1} \text{ and } {C}_n.
   \end{array}
\end{equation}
Note that the second elements of the pairs are obtained by {\it
folding} operation from the first ones, see \cite[Remark 5.4]{St05}.
}
\end{remark}

\section{The orbit structure of the Coxeter transformation}
 \label{orbit_str}
 Let $\mathfrak{g}$ be the simple complex Lie algebra
 of type $A, D$ or $E$, $\tilde{\mathfrak{g}}$ be the affine Kac-Moody
 Lie algebra associated to  $\mathfrak{g}$, and
 $\mathfrak{h} \subseteq \tilde{\mathfrak{h}}$ be, respectively, Cartan
 subalgebras of $\mathfrak{g} \subseteq \tilde{\mathfrak{g}}$.
 Let $\mathfrak{h}^{\vee}$ (resp. ${\tilde{\mathfrak{h}}^{\vee}}$)
 be the dual space to $\mathfrak{h}$ (resp. ${\mathfrak{h}}^{\vee}$)
 and $\alpha_i \in \mathfrak{h}^{\vee}, i = 1,\dots,l$
 be an ordered set of simple positive roots.
 Here, we follow B.~Kostant's description \cite{Kos84}
 of the orbit structure of the Coxeter transformation ${\bf C}$
 on the highest root in the root system of $\mathfrak{g}$.
 We consider a bipartite graph and a bicolored Coxeter
 transformation from \cite[Section 3]{St05}.
 Let $\beta$ be the highest root of ($\mathfrak{h}, \mathfrak{g}$),
 see \cite[Section 4.1]{St05}. Then
 $$
     w_2\beta = \beta \quad  \text{ or } \quad  w_1\beta = \beta.
 $$
In the second case we just swap $w_1$ and $w_2$, i.e., we always
have
 \begin{equation}
   \label{choose_w2}
         w_2\beta = \beta.
 \end{equation}
Between two bicolored Coxeter transformations \cite[(3.1)]{St05} we
select such one that
$$
   {\bf C} = w_2{w}_1.
$$
 Consider, for example, the Dynkin diagram $E_6$. Here,
\begin{equation}
   \label{choose_w1}
    w_1 =  \left (
      \begin{array}{cccccc}
        1 & & &  &  &  \\
        & 1 & &  &  &  \\
        & & 1 &  &  &  \\
        1 & 1 & 0 & -1 & &  \\
        1 & 0 & 1 & & -1 & \\
        1 & 0 & 0 & & & -1 \\
      \end{array}
       \right )
      \begin{array}{c}
         x_0 \\
         x_1 \\
         x_2 \\
         y_1 \\
         y_2 \\
         y_3 \\
      \end{array},
       \hspace{5mm}
   w_2 = \left (
      \begin{array}{cccccc}
        -1 & & & 1 & 1 & 1 \\
        & -1 & & 1 & 0 & 0 \\
        &  & -1 & 0 & 1 & 0 \\
        &  &    & 1 & & \\
        &  &    & & 1 & \\
        &  &    & & & 1 \\
      \end{array}
       \right )
      \begin{array}{c}
         x_0 \\
         x_1 \\
         x_2 \\
         y_1 \\
         y_2 \\
         y_3 \\
      \end{array}
\end{equation}
The vector $z \in \mathfrak{h}^{\vee}$ and the highest root
$\beta$ are:
\begin{equation}
  \label{def_E6_vectors}
      z =
      \begin{array}{ccccc}
        x_1 & y_1 & x_0 & y_2 & x_2 \\
            &     & y_3 &     &
      \end{array},
      \hspace{5mm}
      \beta =
      \begin{array}{ccccc}
        1 & 2 & 3 & 2 & 1 \\
        &     & 2 &     &
      \end{array} , \quad
      \text{ or } \quad
      \beta =
      \left (
      \begin{array}{c}
         3 \\
         1 \\
         1 \\
         2 \\
         2 \\
         2 \\
      \end{array}
      \right )
      \begin{array}{c}
         x_0 \\
         x_1 \\
         x_2 \\
         y_1 \\
         y_2 \\
         y_3 \\
      \end{array}
\end{equation}

Further, following B.Kostant \cite[Th.1.5]{Kos84} consider the
alternating products $\tau^{(n)}$:
\begin{equation}
  \begin{split}
    & \tau^{(1)} = w_1, \\
    & \tau^{(2)} = {\bf C} = w_2w_1, \\
    & \tau^{(3)} = w_1{\bf C} = w_1w_2w_1, \\
    & \dots, \\
    & \tau^{(n)} =
       \begin{cases}
         {\bf C}^k = w_2w_1\dots{w_2}w_1 \text{ for } n = 2k, \\
         w_1{\bf C}^k = w_1w_2w_1\dots{w_2}w_1 \text{ for } n = 2k+1, \\
       \end{cases}
  \end{split}
\end{equation}
and the orbit of the highest root $\beta$ under the action of
$\tau^{(n)}$:

\begin{equation}
   \beta_n = \tau^{(n)}\beta, \quad \text{ where }n = 1,\dots,h-1
\end{equation}
($h$ is the Coxeter number, see \cite[(4.1)]{St05}).

\begin{theorem}{\em\bf (B.~Kostant, \cite[Theorems 1.3, 1.4, 1.5, 1.8]{Kos84})}
  \label{Kostant_vectors_z}

  1)  There exist $z_j \in \tilde{\mathfrak{h}}^{\vee}$,
   where $j = 0,\dots,h$,
   and even integers $2 \le a \le b \le h$
   (see (\ref{Kostant_numbers_a_b}) and Table
   \ref{Kostant_numbers}) such that the generating functions
   $P_{G}(t)$
   (see (\ref{Kostant_gen_func}), (\ref{gen_func_i}) )
   are obtained as follows:
  \begin{equation}
    \label{series_P_G_i}
      [P_{G}(t)]_i =
       \begin{cases}
         \displaystyle\frac{1 + t^h}
            {(1 - t^a)(1 - t^b)} \text{ for } i = 0, \vspace{5mm} \\
         \displaystyle\frac{\sum\limits_{j=0}^{h}z_j{t}^j}
            {(1 - t^a)(1 - t^b)} \text{ for } i = 1,\dots,r. \\
       \end{cases}
  \end{equation}
  For $n = 1,\dots,h-1$, one has
  $z_n \in \mathfrak{h}^{\vee}$ ( not just
  $\tilde{\mathfrak{h}}^{\vee}$).
  The indices $i = 1,\dots,r$ enumerate the vertices of the
  Dynkin diagram and the coordinates of the vectors $z_n$;
  $i = 0$ corresponds to the additional (affine) vertex,
  the one that extends the Dynkin diagram to the extended Dynkin
  diagram. One has $z_0 = z_h = \alpha_0$, where
  $\alpha_0 \in \tilde{\mathfrak{h}}^{\vee}$ is the added simple
  root corresponding to the affine vertex.

  2) The vectors $z_n$ (we call these vectors the {\it assembling vectors})
  are obtained as follows:
\begin{equation}
  \label{assemb_vect}
  z_n = \tau^{(n-1)}\beta - \tau^{(n)}\beta.
\end{equation}

  3) We have
\begin{equation}
  z_g = 2\alpha_{*}, \text{ where } g = \frac{h}{2},
\end{equation}
and $\alpha_{*}$ is the simple root corresponding to the branch
point for diagrams $D_n$, $E_n$ and to the midpoint for the
diagram $A_{2m-1}$. In all these cases $h$ is even, and $g$ is
integer. The diagram $A_{2m}$ has been excluded.

  4) The series of assembling vectors $z_n$ is symmetric:
\begin{equation}
 \label{symmetry_z}
  z_{g+k} = z_{g-k} \text{ for } k = 1,\dots,g.
\end{equation}
\end{theorem}

In the case of the Dynkin diagram $E_6$, the vectors
$\tau^{(n)}\beta$ are given in Table \ref{table_orbit_max_root},
and the assembling vectors $z_n$ are given in Table
\ref{table_vectors_z}. The vector $z_6$ coincides with
$2\alpha_{x_0}$, where $\alpha_{x_0}$ is the simple root
corresponding to the vertex $x_0$, see (\ref{def_E6_vectors}).
From Table \ref{table_vectors_z} we see, that
\begin{equation}
    z_1 = z_{11}, \hspace{5mm}
    z_2 = z_{10}, \hspace{5mm}
    z_3 = z_9 , \hspace{5mm}
    z_4 = z_8, \hspace{5mm}
    z_5 = z_7.
\end{equation}

\begin{table} 
 \centering
 \vspace{2mm}
 \caption{\hspace{3mm}The orbit of the Coxeter transformation on the highest root}
 \renewcommand{\arraystretch}{1.3}
  \begin{tabular} {||c|c|c||}
  \hline \hline
      $\beta$   &
      $\tau^{(1)}\beta = w_1\beta$ &
      $\tau^{(2)}\beta = {\bf C}\beta$ \\
  \hline
      $\begin{array}{ccccc}
        1 & 2 & 3 & 2 & 1 \\
        &     & 2 &     &
      \end{array}$
&     $\begin{array}{ccccc}
        1 & 2 & 3 & 2 & 1 \\
        &     & 1 &     &
      \end{array}$
&     $\begin{array}{ccccc}
        1 & 2 & 2 & 2 & 1 \\
        &     & 1 &     &
      \end{array}$
      \\
  \hline \hline
     $\tau^{(3)}\beta = w_1{\bf C}\beta$ &
     $\tau^{(4)}\beta = {\bf C}^2\beta$ &
     $\tau^{(5)}\beta = w_1{\bf C}^2\beta$
     \\
  \hline
     $\begin{array}{ccccc}
        1 & 1 & 2 & 1 & 1 \\
        &     & 1 &     &
      \end{array}$
&     $\begin{array}{ccccc}
        0 & 1 & 1 & 1 & 0 \\
        &     & 1 &     &
      \end{array}$
&     $\begin{array}{ccccc}
        0 & 0 & 1 & 0 & 0 \\
        &     & 0 &     &
      \end{array}$
      \\
  \hline \hline
      $\tau^{(6)}\beta = {\bf C}^3\beta$ &
      $\tau^{(7)}\beta = w_1{\bf C}^3\beta$ &
      $\tau^{(8)}\beta = {\bf C}^4\beta$ \\
  \hline
      $\begin{array}{ccccc}
        0 & 0 & -1 & 0 & 0 \\
        &     & 0 &     &
      \end{array}$
&     $\begin{array}{ccccc}
        0 & -1 & -1 & -1 & 0 \\
        &     & -1 &     &
      \end{array}$
&      $\begin{array}{ccccc}
        -1 & -1 & -2 & -1 & -1 \\
        &     & -1 &     &
      \end{array}$
      \\
  \hline \hline
      $\tau^{(9)}\beta = w_1{\bf C}^4\beta$  &
      $\tau^{(10)}\beta = {\bf C}^5\beta$ &
      $\tau^{(11)}\beta = w_1{\bf C}^5\beta$ \\
  \hline
      $\begin{array}{ccccc}
        -1 & -2 & -2 & -2 & -1 \\
        &     & -1 &     &
      \end{array}$
&     $\begin{array}{ccccc}
        -1 & -2 & -3 & -2 & -1 \\
        &     & -1 &     &
      \end{array}$
&     $\begin{array}{ccccc}
        -1 & -2 & -3 & -2 & -1 \\
        &     & -2&     &
      \end{array}$ \\
  \hline \hline
\end{tabular}
  \label{table_orbit_max_root}
\end{table}

\begin{table} 
 \centering
 \vspace{2mm}
 \caption{\hspace{3mm}The assembling vectors $z_n = \tau^{(n-1)}\beta -
            \tau^{(n)}\beta$}
 \renewcommand{\arraystretch}{1.3}
  \begin{tabular} {||c|c|c||}
  \hline \hline
      $z_1 = \beta - w_1\beta$   &
      $z_2 = w_1\beta - {\bf C}\beta$ &
      $z_3 = {\bf C}\beta - w_1{\bf C}\beta$ \\
  \hline
      $\begin{array}{ccccc}
        0 & 0 & 0 & 0 & 0 \\
        &     & 1 &     &
      \end{array}$
&     $\begin{array}{ccccc}
        0 & 0 & 1 & 0 & 0 \\
        &     & 0 &     &
      \end{array}$
&     $\begin{array}{ccccc}
        0 & 1 & 0 & 1 & 0 \\
        &     & 0 &     &
      \end{array}$
      \\
  \hline \hline
      $z_4 = w_1{\bf C}\beta - {\bf C}^2\beta$   &
      $z_5 = {\bf C}^2\beta - w_1{\bf C}^2\beta$ &
      $z_6 = w_1{\bf C}^2\beta - {\bf C}^3\beta$ \\
  \hline
      $\begin{array}{ccccc}
        1 & 0 & 1 & 0 & 1 \\
        &     & 0 &     &
      \end{array}$
&     $\begin{array}{ccccc}
        0 & 1 & 0 & 1 & 0 \\
        &     & 1 &     &
      \end{array}$
&     $\begin{array}{ccccc}
        0 & 0 & 2 & 0 & 0 \\
        &     & 0 &     &
      \end{array}$
      \\
  \hline \hline
      $z_7 = {\bf C}^3\beta - w_1{\bf C}^3\beta$   &
      $z_8 = w_1{\bf C}^3\beta - {\bf C}^4\beta$ &
      $z_9 = {\bf C}^4\beta - w_1{\bf C}^4\beta$ \\
  \hline
      $\begin{array}{ccccc}
        0 & 1 & 0 & 1 & 0 \\
        &     & 1 &     &
      \end{array}$
&     $\begin{array}{ccccc}
        1 & 0 & 1 & 0 & 1 \\
        &     & 0 &     &
      \end{array}$
&     $\begin{array}{ccccc}
        0 & 1 & 0 & 1 & 0 \\
        &     & 0 &     &
      \end{array}$
      \\
  \hline \hline
      $z_{10} = w_1{\bf C}^4\beta - {\bf C}^5\beta$   &
      $z_{11} = {\bf C}^5\beta - w_1{\bf C}^5\beta$ & \\
  \hline
      $\begin{array}{ccccc}
        0 & 0 & 1 & 0 & 0 \\
        &     & 0 &     &
      \end{array}$
&     $\begin{array}{ccccc}
        0 & 0 & 0 & 0 & 0 \\
        &     & 1 &     &
      \end{array}$
&     \\
  \hline \hline
\end{tabular}
  \label{table_vectors_z}
\end{table}

Denote by $z(t)_i$ the polynomial $\sum\limits_{j=0}^{h}z_jt^j$ from
(\ref{series_P_G_i}). In the case of $E_6$ we have:
\begin{equation}
 \label{def_zt_2}
  \begin{split}
    & z(t)_{x_0} = t^2 + t^4 + 2t^6 + t^8 + t^{10},  \\
    & z(t)_{x_1} = t^4 + t^8, \\
    & z(t)_{x_2} = t^4 + t^8, \\
    & z(t)_{y_1} = t^3 + t^5 + t^7 + t^9, \\
    & z(t)_{y_2} = t^3 + t^5 + t^7 + t^9, \\
    & z(t)_{y_3} = t + t^5 + t^7 + t^{11}.
  \end{split}
\end{equation}
The Kostant numbers $a,b$ (see Table \ref{Kostant_numbers}) for
$E_6$ are $a = 6$, $b = 8$. From (\ref{series_P_G_i}) and
(\ref{def_zt_2}), we have

\begin{equation}
 \begin{split}
 & [P_G(t)]_{x_0} \hspace{21mm} = \frac{t^2 + t^4 + 2t^6 + t^8 + t^{10}}{(1 - t^6)(1 - t^8)},
    \\ \\
 & [P_G(t)]_{x_1} = [P_G(t)]_{x_2} = \frac{t^4 + t^8}{(1 - t^6)(1 - t^8)},
    \\ \\
 & [P_G(t)]_{y_1} = [P_G(t)]_{y_2} = \frac{t^3 + t^5 + t^7 + t^9}{(1 - t^6)(1 - t^8)},
    \\ \\
 & [P_G(t)]_{y_3} \hspace{21mm} = \frac{t + t^5 + t^7 + t^{11}}{(1 - t^6)(1 - t^8)}.
 \end{split}
\end{equation}
Since
\begin{equation}
     1 - t^6 = \sum\limits_{n=0}^{\infty}t^{6n}, \hspace{5mm}
     1 - t^8 = \sum\limits_{n=0}^{\infty}t^{8n},
\end{equation}
we have
\begin{equation}
 \label{calc_poincare}
 \begin{split}
 & [P_G(t)]_{x_1} = [P_G(t)]_{x_2} =
           \sum\limits_{i,j=0}^{\infty}(t^{6i + 8j + 4} + t^{6i + 8j + 8}),
    \\ \\
 & [P_G(t)]_{y_1} = [P_G(t)]_{y_2} =
          \sum\limits_{i,j=0}^{\infty}(t^{6i + 8j + 3} + t^{6i + 8j + 5} +
                                           t^{6i + 8j + 7} + t^{6i + 8j + 9}),
    \\ \\
 & [P_G(t)]_{y_3} \hspace{21mm} =
          \sum\limits_{i,j=0}^{\infty}(t^{6i + 8j + 1} + t^{6i + 8j + 5}+
                                           t^{6i + 8j + 7} + t^{6i + 8j + 11}),
 \end{split}
\end{equation}

Recall that $m_{\alpha}(n)$, where $\alpha = x_1, x_2, y_1, y_2,
y_3$, are the multiplicities of the indecomposable representations
$\rho_{\alpha}$ of $G$ (considered in the context of the McKay
correspondence, \cite{Kos84}) in the decomposition of $\pi_n|G$
(\ref{decomp_pi_n}). These multiplicities are the coefficients of
the Poincar\'{e} series (\ref{calc_poincare}), see (\ref{def_vn}),
(\ref{Kostant_gen_func}), (\ref{gen_func_i}). For example,
\begin{equation}
  [P_G(t)]_{x_1} = [P_G(t)]_{x_2} =
    t^4 + t^8 + t^{10} + t^{12} + t^{14} + 2t^{16} + t^{18} + 2t^{20} + \dots
\end{equation}
\begin{equation}
  m_{x_1} = m_{x_2} =
  \begin{cases}
     0 ~\text{ for } n = 1,2,3,5,6 \text{ and } n = 2k+1, k \geq 3, \\
     1 ~\text{ for } n = 4,8,10,12,14,18, \dots \\
     2 ~\text{ for } n = 16,20, \dots \\
     \dots
  \end{cases}
\end{equation}
In particular, the representations $\rho_{x_1}(n)$ and
$\rho_{x_2}(n)$ do not enter in the decomposition of $\pi_n$ of
$SU(2)$ (see \S\ref{generating_fun}) for all odd $n$.

  In \cite{Kos04}, concerning the importance of the polynomials
  $z(t)_i$, B.~Kostant points out:
  {\it ``Unrelated to the Coxeter element, the polynomials $z(t)_i$
are also determined in Springer, \cite{Sp87}. They also appear in
another context in Lusztig, \cite{Lus83} and \cite{Lus99}. Recently,
in a beautiful result, Rossmann, \cite{Ros04}, relates the character
of $\gamma_i$ to the polynomial $z(t)_i$.''}

\section{McKay's observation relating the Poincar\'{e} series}
  \label{McKay_observ}
  In this section we prove McKay's observation \cite{McK99}
  relating the Poincar\'{e} series, or rather Molien-Poincar\'{e} series.
  In our context these series are the Kostant generating functions $P_G(t)$
  corresponding the indecomposable representations of group $G$:
 \begin{equation}
   \label{Moilen_zt}
    [P_G(t)]_i = \frac{z(t)_i}{(1-t^a)(1-t^b)},
 \end{equation}
 see (\ref{series_P_G_i}), (\ref{def_zt_2}).

 \begin{theorem}{\em (McKay's observation \cite[(*)]{McK99})}
 For diagrams $G = D_n, E_n$ and $A_{2m-1}$,
 the Kostant generating functions $[P_G(t)]_i$ are related as
 follows:
 \begin{equation}
   \label{2_McKay}
    (t+t^{-1})[P_G(t)]_i = \sum\limits_{i \leftarrow j}[P_G(t)]_j,
 \end{equation}
where $j$ runs over all vertices adjacent to $i$, and $[P_G(t)]_0$
related to the affine vertex $\alpha_0$ occurs in the right side
only: $i = 1,2,\dots,r$.
 \end{theorem}

 By (\ref{Moilen_zt}), McKay's observation (\ref{2_McKay}) is
equivalent to the following one:
 \begin{equation}
   \label{2_McKay_1}
    (t+t^{-1})z(t)_i = \sum\limits_{i \leftarrow j}z(t)_j, \text{
    where } i = 1,\dots,r.
 \end{equation}
  So, we will prove (\ref{2_McKay_1}).

  The {\it adjacency matrix} $\mathcal{A}$ for types $ADE$ is the
  matrix containing non-diagonal entries $a_{ij}$ if and
  only if the vertices $i$ and $j$ are connected by an edge, and
  then $a_{ij} = 1$, and all diagonal entries $a_{ii}$ vanish:
\begin{equation}
   \label{adj_matr_A}
     \mathcal{A} = \left(
          \begin{array}{cc}
             0 & -2D \\
             -2D^t & 0
          \end{array}
         \right).
 \end{equation}

\begin{figure}[h]
\centering
\includegraphics{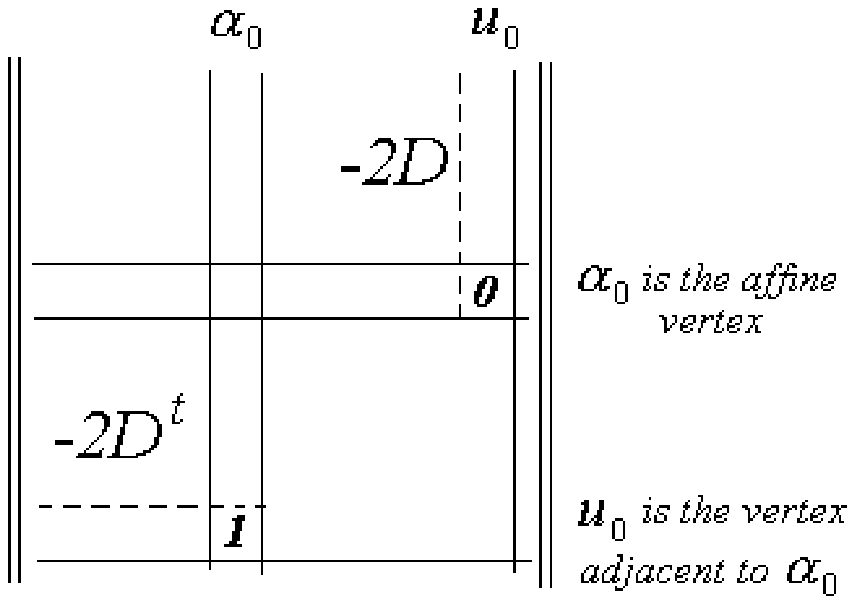} \caption{\hspace{3mm}
The semi-affine adjacent matrix $\mathcal{A}^{\gamma}$ }
\label{semi_affine}
\end{figure}

  Let $\alpha_0$ be the affine vertex of the graph
  $\varGamma$, and $u_0$ be a vertex adjacent to $\alpha_0$.
  Extend the adjacency matrix $\mathcal{A}$ to the {\it semi-affine adjacency
  matrix} $\mathcal{A}^{\gamma}$ (in the style to the McKay definition of
  the semi-affine graph in \cite{McK99}) by adding
  a row and a column corresponding to the affine vertex
  $\alpha_0$ as follows: $0$ is set in the $(u_0, \alpha_0)$th
  slot and $1$ is set in the$(\alpha_0, u_0)$th slot, all remaining
  places in the $\alpha_0$th row  and the $\alpha_0$th column are $0$,
  see Fig. \ref{semi_affine}. Note, that for the
  $A_n$  case, we set $1$ in two places: $(\alpha_0, u_0)$ and $(\alpha_0, u'_0)$
  corresponding to vertices $u_0$ and $u'_0$  adjacent to $\alpha$.

Using the semi-affine adjacency matrix $\mathcal{A}^{\gamma}$ we
write McKay's observation (\ref{2_McKay_1}) in the matrix form:
 \begin{equation}
   \label{2_McKay_2}
     (t + t^{-1})z(t)_i = (\mathcal{A}^{\gamma}z(t))_i, \quad \text{ where }
     z(t) = \{z(t)_0,\dots,z(t)_r \} \text{ and } i = 1,\dots,r.
 \end{equation}

  To prove (\ref{2_McKay_2}), we consider the action of
  the adjacency matrix $\mathcal{A}$
  and the semi-affine adjacency matrix $\mathcal{A}^{\gamma}$
  related to the extended Dynkin diagram of types $ADE$ on
  assembling vectors $z_n$ (\ref{assemb_vect}).

\begin{proposition}
  1) For the vectors $z_i \in {\mathfrak{h}}^{\vee}$ from (\ref{assemb_vect}), we have
  \begin{equation}
    \label{adj_1}
    \begin{split}
      & (a) \quad  \mathcal{A}z_i =
          z_{i-1} + z_{i+1} \quad \text{ for } \quad 1 < i < h-1, \\
      & (b) \quad  \mathcal{A}z_1 =
          z_2 \quad \text { and } \quad \mathcal{A}z_{h-1} = z_{h-2}, \\
    \end{split}
  \end{equation}

  2) Consider the same vectors $z_i$ as vectors from
  ${\tilde{\mathfrak{h}}}^{\vee}$, so we just add a
  zero coordinate to the affine vertex $\alpha_0$. We have
  \begin{equation}
    \label{adj_1_ext}
    \begin{split}
      & (a) \quad \mathcal{A}^{\gamma}z_i = z_{i-1} + z_{i+1}
        \quad \text{ for } \quad 1 < i < h-1, \\
      & (b) \quad \mathcal{A}^{\gamma}z_1 = z_2
        \quad \text { and } \quad \mathcal{A}^{\gamma}z_{h-1} = z_{h-2},\\
      & (c) \quad \mathcal{A}^{\gamma}z_0 = z_1 \quad \text { and }
        \quad  \mathcal{A}^{\gamma}z_h = z_{h-1}.
    \end{split}
  \end{equation}
\end{proposition}

\PerfProof

  1a) Let us prove (\ref{adj_1} a).
  According to (\ref{assemb_vect}) we have
  \begin{equation}
     z_{2n} = w_{1}{\bf C}^{n-1}\beta - {\bf C}^n\beta =
      (1 - w_2)w_1{\bf C}^{n-1}\beta,
  \end{equation}
  and
  \begin{equation}
     z_{2n+1} = {\bf C}^n\beta - w_1{\bf C}^n\beta =
      (1 - w_1){\bf C}^{n}\beta.
  \end{equation}
 Thus, for $i = 2n$, eq. (\ref{adj_1} a) is equivalent to
  \begin{equation}
    \label{adj_even_1}
    \mathcal{A}(1 - w_2)w_1{\bf C}^{n-1}\beta =
        (1 - w_1){\bf C}^{n-1}\beta + (1 - w_1){\bf C}^{n}\beta,
  \end{equation}
 and for $i = 2n + 1$, eq. (\ref{adj_1} a) is equivalent to
  \begin{equation}
    \label{adj_odd_1}
    \mathcal{A}(1 - w_1){\bf C}^{n}\beta =
        (1 - w_2)w_1{\bf C}^{n-1}\beta +
         (1 - w_2)w_1{\bf C}^{n}\beta.
  \end{equation}
  To prove relations (\ref{adj_even_1}) and (\ref{adj_odd_1}), it
  suffices to show that
  \begin{equation}
    \label{adj_even_2}
    \mathcal{A}(1 - w_2)w_1 = (1 - w_1) + (1 - w_1){\bf C} =
      (1 - w_1)(1 + {\bf C}),
  \end{equation}
  and
  \begin{equation}
    \label{adj_odd_2}
    \mathcal{A}(1 - w_1){\bf C} = (1 - w_2)w_1 + (1 - w_2)w_1{\bf C} =
      (1 - w_2)w_1(1 + {\bf C}).
  \end{equation}
 In (\ref{choose_w2}), (\ref{choose_w1}), $w_1$ and $w_2$ are chosen as
\begin{equation}
 w_1 = \left(
  \begin{array}{cc}
     I    &  0 \\
    -2D^t & -I
  \end{array}
 \right), \qquad
 w_2 = \left(
  \begin{array}{cc}
     -I    &  -2D \\
      0  &  I
  \end{array}
 \right),
\end{equation}
So, by (\ref{adj_matr_A}) we have
\begin{equation}
  \label{6_matrix}
 \begin{array}{cc}
 & 1 - w_1 = \left(
  \begin{array}{cc}
     0    &  0 \\
    2D^t & 2I
  \end{array}
 \right), \quad
 1 - w_2 = \left(
  \begin{array}{cc}
     2I  &  2D \\
      0  &  0
  \end{array}
 \right),   \vspace{3mm} \\
 & (1 - w_2)w_1 = \left(
  \begin{array}{cc}
  2I - 4DD^t  &  -2D \\
    0         &  0
  \end{array}
   \right), \quad
 {\bf C} = \left(
  \begin{array}{cc}
     4DD^t - I  &  2D \\
      -2D^t  &  -I
  \end{array}
 \right), \vspace{3mm} \\
 & \mathcal{A}(1 - w_2)w_1 = \left(
  \begin{array}{cc}
    0                &   0 \\
  - 4D^t + 8D^tDD^t  &  4D^tD
  \end{array}
   \right), \vspace{3mm} \\
 & 1 + {\bf C} = \left(
  \begin{array}{cc}
     4DD^t   &  2D \\
      -2D^t  &  0
  \end{array}
 \right).
  \end{array}
\end{equation}
By (\ref{6_matrix}) we have
\begin{equation}
 (1 - w_1)(1 + {\bf C}) = \left(
  \begin{array}{cc}
    0                &   0 \\
  - 4D^t + 8D^tDD^t  &  4D^tD
  \end{array}
 \right),
\end{equation}
 and (\ref{adj_even_2}) is true.
Further, we have
\begin{equation}
  \label{more_matr_1}
 \begin{array}{cc}
 & (1 - w_1){\bf C} = \left(
  \begin{array}{cc}
     0    &  0 \\
    8D^tDD^t - 6D^t & 4D^tD - 2I
  \end{array}
 \right), \vspace{3mm} \\
 & \mathcal{A}(1 - w_1){\bf C} = \left(
  \begin{array}{cc}
      -16DD^tDD^t +12DD^t &  -8DD^tD + 4D \\
      0                  &  0
  \end{array}
 \right),   \vspace{3mm} \\
 \end{array}
\end{equation}
By (\ref{6_matrix}) we obtain
 \begin{equation}
 (1 - w_2)w_1(1 + {\bf C}) = \left(
  \begin{array}{cc}
      -16DD^tDD^t +12DD^t &  -8DD^tD + 4D \\
      0                  &  0
  \end{array}
 \right),   \vspace{3mm} \\
\end{equation}
and (\ref{adj_odd_1}) is also true.

  1b) Let us move on to (\ref{adj_1} b). This is equivalent to
\begin{equation}
  \label{A_z1_z2}
  \begin{split}
   & \mathcal{A}(\beta - w_1\beta) =
      w_1\beta - {\bf C}\beta, \quad  \text{or} \\
   & \mathcal{A}(1 - w_1)\beta = (1 - w_2)w_1\beta.
  \end{split}
\end{equation}
By (\ref{adj_matr_A}), (\ref{6_matrix}) eq. (\ref{A_z1_z2}) is
equivalent to
\begin{equation}
  \label{A_z1_z2_2}
  \left (
  \begin{array}{cc}
    -4DD^t & -4D \\
    0      &   0
  \end{array}
  \right ) \left (
  \begin{array}{c}
    x \\
    y
  \end{array}
  \right ) =
  \left (
  \begin{array}{cc}
    2I - 4DD^t & -2D \\
    0      &   0
  \end{array}
  \right ) \left (
  \begin{array}{c}
    x \\
    y
  \end{array}
  \right ),
\end{equation}
where
 \begin{equation}
  \label{partition_xy}
 \beta = \left (
  \begin{array}{c}
    x \\
    y
  \end{array}
  \right )
\end{equation}
is given in two-component form corresponding to a bipartite graph.
Eq. (\ref{A_z1_z2_2}) is equivalent to
\begin{equation}
  \label{adjac_rel}
  -2Dy = 2x .
\end{equation}
 Since the matrix $-2D$ contains a $1$ at the $(i,j)$th slot if and
only if the vertices $i$ and $j$ are connected, eq.
(\ref{adjac_rel}) follows from the well-known fact formulated in
 Remark \ref{highest_vec}(b).

\begin{remark}[On the highest root and imaginary roots]
  \label{highest_vec}
  {\rm
  (a) The highest root $\beta$ for types $ADE$ coincides with the
  minimal positive {\it imaginary root} (of the corresponding extended
  Dynkin diagram) without affine coordinate $\alpha_0$ ,
  see \cite[Section 2.4]{St05}. The coordinates of the
  imaginary vector $\delta$ are given on \cite[Fig. 2.2]{St05}.

  (b) For any vertex $x_i$ of the the highest root $\beta$
   (except the vertex $u_0$ adjacent to the
   affine vertex $\alpha_0$), the sum of coordinates
   of the adjacent vertices $y_j$ coincides
   with the doubled coordinate of $\alpha_{x_i}$:

   \begin{equation}
     \label{doubled_coord}
     \sum\limits_{y_i \rightarrow x_j}{\alpha_{y_j}} =
     2\alpha_{x_i}.
   \end{equation}

  (c) For any vertex $x_i$ of the the imaginary root $\delta$,
   the sum of coordinates of the adjacent vertices $y_j$ coincides
   with the doubled coordinate of $\alpha_{x_i}$ as above in
   (\ref{doubled_coord}).
  }
\end{remark}

  Since we can choose partition (\ref{partition_xy}) such that
  $\alpha_0$ belongs to subset $y$, we obtain (\ref{adjac_rel}).
  The second relation of (\ref{adj_1} b) follows from the symmetry
  of assembling vectors, see (\ref{symmetry_z})
  of the Kostant theorem (Theorem \ref{Kostant_vectors_z}).
  \qed

 2) Relations (a), (b) of (\ref{adj_1_ext})
    follow from the corresponding relations
    in (\ref{adj_1}) since the addition of a ``1'' to the  $(\alpha_0,
    u_0)$th slot of the matrix $\mathcal{A}^{\gamma}$ (\ref{semi_affine})
    is neutralized by the affine coordinate 0 of the vectors $z_i$,
    where $i = 1,\dots,h-1$.

    Let us prove (c) of (\ref{adj_1_ext}). First,
 \begin{equation}
   \begin{split}
    & \mathcal{A}^{\gamma}z_0 = \alpha_{u_0}, \quad \text{ for } D_n, E_n, \\
    & \mathcal{A}^{\gamma}z_0 =
       \alpha_{u_0} + \alpha_{u'_0} , \quad \text{ for } A_{2m-1},
   \end{split}
 \end{equation}
  where $\alpha_{u_0}$ (resp. $\alpha_{u'_0}$)
  is the simple root with a ``1'' in the
  $u_0$th (resp. $u'_0$th) position.
    Thus, by (\ref{6_matrix}) (c) is equivalent to
 \begin{equation}
  \begin{split}
    & \alpha_{u_0} = (1 - w_{1})\beta =
    \left (
    \begin{array}{cc}
       0 \\
       2D^t{x} + 2y
     \end{array} \right ), \quad \text{ for } D_n, E_n,
     \\
    & \alpha_{u_0} +  \alpha_{u'_0} = (1 - w_{1})\beta =
    \left (
    \begin{array}{cc}
       0 \\
       2D^t{x} + 2y
     \end{array} \right ), \quad \text{ for } A_{2m-1},
   \end{split}
 \end{equation}

Again, by Remark \ref{highest_vec} (b), we have $2D^t{x} + 2y = 0$
for all coordinates excepting coordinate $u_0, u'_0$. For
coordinate $u_0$ (resp. $u'_0$), by Remark \ref{highest_vec} (c),
we have $(\alpha_{u_0} - 2D^t{x})_{u_0} = 2y_{u_0}$
 (resp. $(\alpha_{u'_0} - 2D^t{x})_{u'_0} = 2y_{u'_0}$ for $A_{2m-1}$).
 \qedsymbol

 {\it Proof of (\ref{2_McKay_2}) }.  Since,
$$
   z(t) = \sum\limits_{j=0}^{h}z_j{t}^j, \quad
   z(t)_i = (\sum\limits_{j=0}^{h}z_j{t}^j)_i, \text{ where } i =
   1,\dots,r,
$$
by (\ref{adj_1_ext}) we have
\begin{equation*}
  \begin{split}
    \mathcal{A}^{\gamma}z(t) & =
       \sum\limits_{j=0}^{h}\mathcal{A}^{\gamma}z_j{t}^j =\\
    z_1 + z_2{t} + &(z_1 + z_3)t^2 + \dots + \\
      & (z_{h-3} + z_{h-1})t^{h-2} + z_{h-2}t^{h-1} + z_{h-1}t^h = \\
    (z_1 + z_2{t} + & z_3{t}^2 + \dots + z_{h-1}t^{h-2}) + \\
      & (z_1{t}^2 + \dots  + z_{h-2}t^{h-1} + z_{h-1}t^h) = \\
    t^{-1}(z_1{t} + & z_2{t}^2  + z_3{t}^3 + \dots + z_{h-1}t^{h-1}) + \\
      & t(z_1{t} + \dots +  z_{h-2}t^{h-2} + z_{h-1}t^{h-1}) = \\
    & (t + t^{-1})(z_1{t} + z_2{t}^2 + z_3{t}^3 + \dots + z_{h-1}t^{h-1}) = \\
    & (t + t^{-1})(z(t) - z_0 - z_h{t}^h).
  \end{split}
\end{equation*}
 Since $z_0 = z_h$ we have
\begin{equation}
   \mathcal{A}^{\gamma}z(t) = (t + t^{-1})z(t) - (t + t^{-1})(1+t^h)z_0.
\end{equation}
Coordinates $(z_0)_i = (z_h)_i$ are zeros for $i =1,\dots,r$, and
\begin{equation}
  \label{A_gamma_rel}
   (\mathcal{A}^{\gamma}z(t))_i = (t + t^{-1})z(t)_i, \quad i = 1,\dots,r.
\end{equation}
For the coordinate $i = 0$, corresponding to affine vertex
$\alpha_0$, by definition of $\mathcal{A}^{\gamma}$ (see Fig.
\ref{semi_affine}) we have $(\mathcal{A}^{\gamma}z(t))_0 = 0$, and
 $z(t)_0 = (1 + t^h)z_0$.
 \qedsymbol

  Let us check McKay's observation for the Kostant generating functions
  for the case of $E_6$. According to (\ref{def_zt_2})
  we should get the following relations:

\begin{equation}
  \begin{split}
  & 1) \text{ For } x_0: \quad
        (t+t^{-1})z(t)_{x_0} = z(t)_{y_1} + z(t)_{y_2} + z(t)_{y_3}, \\
  & 2) \text{ For } y_1: \quad (t+t^{-1})z(t)_{y_1} = z(t)_{x_0} + z(t)_{x_1}, \\
  & 3) \text{ For } y_2: \quad (t+t^{-1})z(t)_{y_2} = z(t)_{x_0} + z(t)_{x_2}, \\
  & 4) \text{ For } x_1: \quad (t+t^{-1})z(t)_{x_1} = z(t)_{y_1}, \\
  & 5) \text{ For } x_2: \quad (t+t^{-1})z(t)_{x_2} = z(t)_{y_2}, \\
  & 6) \text{ For } y_3: \quad (t+t^{-1})z(t)_{y_3} = z(t)_{x_0} + z(t)_{\alpha_0}. \\
  \end{split}
\end{equation}
 1) For $x_0$, we have
\begin{equation}
  \begin{split}
    & (t+t^{-1})(t^2 + t^4 + 2t^6 + t^8 + t^{10}) = \\
    & 2(t^3 + t^5 + t^7 + t^9) + (t + t^5 + t^7 + t^{11}), \\
    & \text{ or }\\
    & (t^3 + t^5 + 2t^7 + t^9 + t^{11}) +
      (t + t^3 + 2t^5 + t^7 + t^9) = \\
    & t + 2t^3 + 3t^5 +3t^7 + 2t^9 + t^{11}.
  \end{split}
\end{equation}
 2) For $y_1$ (the same for $y_2$), we have
\begin{equation}
  \begin{split}
    & (t+t^{-1})(t^3 + t^5 + t^7 + t^9) = \\
    &   (t^2 + t^4 + 2t^6 + t^8 + t^{10}) + (t^4 + t^8), \\
    & \text{ or }\\
    & (t^4 + t^6 + t^8 + t^{10}) + (t^2 + t^4 + t^6 + t^{8}) = \\
    &  (t^2 + 2t^4 + 2t^6 + 2t^8 + t^{10}).
  \end{split}
\end{equation}
 3) For $x_1$ (the same for $x_2$), we have
\begin{equation}
    (t+t^{-1})(t^4 + t^8) = t^3 + t^5 + t^7 + t^9.
\end{equation}
 4) For $y_3$, we have
\begin{equation}
  \begin{split}
    & (t+t^{-1})(t + t^5 + t^7 + t^{11}) = \\
    & (t^2 + t^4 + 2t^6 + t^8 + t^{10}) + (1 + t^{12}), \\
    & \text{ or }\\
    & (t^2 + t^6 + t^8 + t^{12}) + (1 + t^4 + t^6 + t^{10}) = \\
    & (t^2 + t^4 + 2t^6 + t^8 + t^{10}) + (1 + t^{12}).
  \end{split}
\end{equation}